\documentclass{amsart}
\usepackage{amsthm}
\usepackage{amsmath}
\usepackage{amstext}
\usepackage{amsfonts}
\usepackage{latexsym}
\usepackage{amscd}
\usepackage{xypic}
\newcommand{\la}{\ensuremath{\rightarrow}}
\theoremstyle{plain}
\newtheorem{theorem}{Theorem}[section]
\newtheorem{lemma}[theorem]{Lemma}
\newtheorem{proposition}[theorem]{Proposition}
\newtheorem{corollary}[theorem]{Corollary}

\theoremstyle{definition}
\newtheorem{definition}[theorem]{Definition}

\newtheorem{example}[theorem]{Example}
\newtheorem{remark}[theorem]{Remark}
\newtheorem{question}[theorem]{Question}
\numberwithin{equation}{theorem}

\begin{document}

\title[Smoothings of normal crossing Fano schemes.]{Smoothings of Fano schemes with normal crossing singularities of dimension at most three.}
\author{Nikolaos Tziolas}
\address{Department of Mathematics, University of Cyprus, P.O. Box 20537, Nicosia, 1678, Cyprus}
\email{tziolas@ucy.ac.cy}
\thanks{This paper was partially written during the author's stay at the Max Planck Institute for Mathematics in Bonn, from April to August 2009.}

\subjclass[2000]{Primary 14D15, 14D06, 14J45.}



\begin{abstract}
In this paper we study the deformation theory of a Fano variety $X$ with normal crossing singularities. We obtain a formula for $T^1(X)$ in a suitable log resolution of $X$ and we obtain explicit criteria for the existence of smoothings of $X$.
\end{abstract}

\maketitle

\section{Introduction}
In this paper we study the deformation theory of a Fano variety of dimension at most three with normal crossing singularities. In particular we investigate when such a variety is smoothable. By this we mean that there is a flat  projective morphism $ f \colon \mathcal{X} \la \Delta$, where $\Delta$ is the spectrum of a discrete valuation ring $(R,m_R)$, such that $\mathcal{X} \otimes_R (R/m_R) \cong X$ and $\mathcal{X} \otimes_R K(R)$ is smooth over $K(R)$. Moreover, we study when such a smoothing exists with smooth total space $\mathcal{X}$. In this case we say that $X$ is totally smoothable.

Normal crossing singularities appear quite naturally in any degeneration problem. Let $f \colon \mathcal{X} \la C$ be a flat projective morphism from a variety $\mathcal{X}$ to a curve $C$. Then, according to Mumford's semistable reduction theorem~\cite{KKMS73},  after a finite base change and a birational modification the family can be brought to standard form $f^{\prime} \colon \mathcal{X}^{\prime} \la C^{\prime}$, where $\mathcal{X}^{\prime}$ is smooth and the special fibers are simple normal crossing varieties.

Smoothings of Fano varieties play a fundamental role in higher dimensional birational geometry. The outcome of the minimal model program starting with a smooth $n$-dimensional projective variety $X$, is a $\mathbb{Q}$-factorial terminal projective variety $Y$ such that either $K_Y$ is nef, or $Y$ has a Mori fiber space structure. This means that there is a projective morphism $f \colon Y \la Z$ such that $-K_Y$ $f$-ample, $Z$ is normal and $\dim Z \leq \dim X -1$. Suppose that the second case happens and $\dim Z=1$. Let $z \in Z$ and $Y_z= f^{-1}(z)$. Then $Y_z$ is a Fano variety of dimension $n-1$ and $Y$ is a smoothing of $Y_z$. The singularities of the special fibers are difficult to describe but normal crossing singularities naturally occur and are the simplest possible non-normal singularities. 

Moreover, the study of smoothings $f \colon \mathcal{X} \la \Delta$ such that $\mathcal{X}$ is smooth, $-K_{\mathcal{X}}$ is $f$-ample and the special fiber is a simple normal crossing divisor, has a central role in the classification of smooth Fano varieties~\cite{Fu90}. In dimension two T. Fujita~\cite{Fu90} has described all the possible degenerations of smooth Del Pezzo surfaces to simple normal crossing Del Pezzo surfaces and Y. Kachi showed that all these actually occur~\cite{Kac07}. As far as we know this problem is completely open in higher dimensions. 

It is therefore of interest to study which Fano varieties with normal crossing singularities are smoothable. Moreover, in view of the second problem mentioned above, the case of simple normal crossing Fano varieties deserves special attention too.

This paper is organized as follows.

In section 3 we describe the tangent space of the versal deformation space $\mathrm{Def}(X)$ of a variety $X$ with normal crossing singularities. In particular, we obtain formulas for the sheaf of first order deformations $T^1(X)$ of $X$. Since $X$ has normal crossing singularities, $T^1(X)$ is a line bundle on its singular locus $Z$. The main result is the following.

\begin{theorem}
Let $X$ be a scheme with normal crossing singularities with $\dim X \leq 3$ and let $Z$ be its singular locus. Let $\phi \colon \tilde{X} \la X$ be a birational morphism such that
\begin{enumerate}
\item $(\tilde{X},1/2\tilde{Z})$, $\tilde{Z}$ are terminal;
\item $K_{\tilde{X}}+1/2\tilde{Z}$ and $K_{\tilde{Z}}$ are $\phi$-nef.
\end{enumerate}
where $\tilde{Z} \subset \tilde{X}$ be the reduced divisorial part of $\phi^{-1}(Z)$ that dominates $Z$ (such spaces do exist). Then
\[
\phi^{\ast}T^1(X)=\mathcal{O}_{\tilde{Z}}(\tilde{Z})^{\otimes 2} \otimes \mathcal{O}_{\tilde{X}}(3E)\otimes \mathcal{O}_{\tilde{Z}}
\]
where $E\subset \tilde{X}$ is the reduced $\phi$-exceptional divisor that dominates the set of singular points of $X$ of multiplicity at least three.
\end{theorem}
In dimension 2, $\tilde{Z} \subset \tilde{X}$ is just the minimal log resolution of $Z \subset X$ and the previous result is a special case of~\cite[Theorem 3.1]{Tzi09}. To prove it we use the explicit minimal model program in dimension three and hence the difficulty to extend it to higher dimensions. In Theorem~\ref{2} we give a formula for $T^1(X)$, where $X$ is of any dimension, in a suitable log resolution $\tilde{Z} \subset \tilde{X}$ of $Z \subset X$. However, the disadvantage of this is that $\tilde{X}$ is not determined by any numerical data.

In the remaining sections of the paper we concentrate on the study of Fano varieties with normal crossing singularities.

In section 4 we study the obstruction spaces for deforming a Fano variety $X$ with normal crossing singularities. It is well known that $H^2(T_X)$ and $H^1(T^1(X))$ are obstruction spaces for deformations of $X$. If $X$ is simple normal crossing, which means that $X$ has smooth irreducible components, then its obstruction theory is deeply clarified by the work of Friedman~\cite{Fr83}. However, in the general case when $X$ is not necessarily reducible, Friedman's theory does not directly apply. In Theorem~\ref{ob4} we show that if $X$ is a Fano variety with normal crossing singularities and $\dim X \leq 3$, then $H^2(T_X)=0$. However, If $X$ has at worst double points then its obstruction theory is much easier to describe. In Theorem~\ref{db-ob} we show that $H^2(T_X)=0$ in all dimensions. As far as the other obstruction space $H^1(T^1(X))$ is concerned, in general it does not vanish. This is the case in example~\ref{ex2}. However, $T^1(X)$ is a line bundle on the singular locus $Z$ of $X$ and in order to obtain smoothings of $X$ one has to impose some positivity conditions on $T^1(X)$ that will force it to vanish.

In sections 5 and 6 we apply the deformation theory developed in the previous sections to obtain criteria for the existence of a smoothing of a Fano variety $X$ with normal crossing singularities. Special attention is given to the case when $X$ is a simple normal crossing Fano variety. If this is the case, Proposition~\ref{sm1} shows that if $X$ has a smoothing with smooth total space then $T^1(X) \cong \mathcal{O}_Z$. This is equivalent to say that $X$ satisfies the $d$-semistability condition defined by Friedman~\cite{Fr83}. Then $X$ admits a logarithmic structure~\cite{KawNa94},~\cite{Kat96} and so we can use the theory of logarithmic deformations developed by K. Kato~\cite{Ka88} and Kawamata-Namikawa~\cite{KawNa94} to study its deformation theory. The main result of sections 5 and 6 is the following.

\begin{theorem}
Let $X$ be a Fano variety with normal crossing singularities. Then
\begin{enumerate}
\item Suppose that $\dim X \leq 3$. Then if $T^1(X)$ is finitely generated by global sections and $H^1(T^1(X))=0$, $X$ is smoothable.
\item Suppose that $X$ is of any dimension and it has only double points. Then if $T^1(X)$ is finitely generated by its global sections, $X$ is smoothable.
\item Suppose that $X$ is simple normal crossing. Then $X$ is totally smoothable if and only if $T^1(X) \cong \mathcal{O}_Z$, where $Z$ is the singular locus of $X$.
\end{enumerate}
\end{theorem}
We do not know if the condition that $T^1(X)$ is finitely generated by its global sections is a necessary condition too for $X$ to be smoothable. In all the cases of the previous theorem it implies that $\mathrm{Def}(X)$ is smooth. If it is true that $\mathrm{Def}(X)$ is smooth for any $X$, then $X$ smoothable implies that $T^1(X)$ is finitely generated by its global sections too.
 
The reason that we cannot extend the previous theorem in all dimensions, is our inability to show that $H^2(T_X)=0$ for a general $X$ of dimension at least 4 having singularities of multiplicity at least three. Proposition~\ref{ob3} contains the main technical result that allows us to show $H^2(T_X)=0$ in the cases of the previous theorem. We believe that it should be possible to generalize it to higher dimensions, but at the moment we have some technical difficulties to do so.

Finally in section 7 we give an example of a smoothable and one of a non-smoothable Fano threefold.

\section{Preliminaries.}

All schemes in this paper are defined over the field of complex numbers $\mathbb{C}$.

A reduced scheme $X$ of finite type over $\mathbb{C}$ is called a normal crossing variety of dimension $n$ if for any $P \in X$, $\hat{\mathcal{O}}_{X,P} \cong \mathbb{C}[[x_0,\dots,x_n]]/(x_0\cdots x_r)$, for some $r=r(P)$, where $\hat{\mathcal{O}}_{X,P}$ is the completion of the local ring of $X$ at $P$ at its maximal ideal. In addition, if $X$ has smooth irreducible components then it is called a simple normal crossing variety.

A reduced projective scheme $X$ with normal crossing singularities is called a Fano variety if and only if $\omega_X^{-1}$ is an ample invertible sheaf on $X$.

For any scheme $X$ we denote by $T^1(X)$ the sheaf of first order deformations of $X$~\cite{Sch68}. If $X$ is reduced then $T^1(X)=\mathcal{E}xt^1_X(\Omega_X,\mathcal{O}_X)$.

Let $X$ be a simple normal crossing variety, let $X_j$ ($1 \leq j \leq k$) its irreducible components and $I_{X_j}$ the ideal sheaves of $X_j$ in $X$. Let $Z$ be its singular locus. Then $X$ is called $d$-semistable~\cite{Fr83} if and only if
\[
(I_{X_1}/I_{X_1}I_Z) \otimes \cdots  \otimes (I_{X_k}/I_{X_k}I_Z) \cong  \mathcal{O}_Z
\]
It is well known that any $d$-semistable simple normal crossing variety  admits a logarithmic structure~\cite{KawNa94},~\cite{Kat96} and so 
in order to investigate the existence of smoothings of such a variety we can use the theory of logarithmic deformations. For definitions and basic properties
we refer the reader to~\cite{KawNa94} for the complex analytic case and~\cite{Ka88},~\cite{Kat96} for the algebraic case. For the convenience of the reader we state the main result that will be used in this paper.

\begin{theorem}[~\cite{KawNa94}]
Let $X$ be a $d$-semistable normal crossing variety and let $(X,\ \mathcal{U})$ be a logarithmic structure on it. Let $m$ be the number of irreducible components of $X$ and $\Lambda_m=\mathbb{C}[[x_1,\ldots,x_m]]$. Let 
\[
0 \la J \la B \la A \la 0
\]
be an extension of local Artin $\Lambda_m$-algebras and let $(X_A,\mathcal{U}_A)$ be a logarithmic deformation of $(X,\mathcal{U})$ over $A$ with structure map $f \colon X_A \la \mathrm{Spec} A $. Then there exists an intrinsically defined locally free $\mathcal{O}_{X_A}$-module $\Omega_{X_A/A}(\mathrm{log})$ of rank $n=\dim X$ such that the obstruction to the existence of a lifting of  $(X_A,\mathcal{U}_A)$ over $B$ is in $H^2(X_A,T_{X_A/A}(\mathrm{log}) \otimes_A J)$, where $T_{X_A/A}(\mathrm{log})=\mathcal{H}om_{X_A}(\Omega_{X_A/A}(\mathrm{log}), \mathcal{O}_{X_A})$. In particular, if $H^2(X,T_{X/\mathbb{C}}(\mathrm{log}))=0$, then $X$ is smoothable by a flat deformation.
\end{theorem}
The previous theorem was proven by Y. Kawamata and Y. Namikawa for complex analytic spaces. However, the results of K. Kato~\cite{Ka88} and F. Kato~\cite{Kat96} show that it also holds in the algebraic category too.

We will repeatedly make use of the Akizuki-Kodaira-Nakano vanishing theorem and its logarithmic version, which we state next.

\begin{theorem}[Akizuki-Kodaira-Nakano~\cite{AN54},~\cite{EV92}]\label{KN}
Let $X$ be a smooth variety and $\mathcal{L}$ an ample invertible sheaf on $X$. Then
\item
\[
H^b(\Omega_X^a\otimes \mathcal{L}^{-1})=0
\]
for all $a,b$ such that $a+b < \dim X$.
\item Moreover, if $D$ is a reduced simple normal crossings divisor of $X$, then
\[
H^b(\Omega_X^a(\mathrm{log}(D))\otimes \mathcal{L}^{-1})=0
\]
for all $a,b$ such that $a+b < \dim X$.
\end{theorem} 
\section{Description of $T^1$}
In this section we describe the sheaf $T^1(X)$ of first order deformations of a scheme $X$ with normal crossing singularities.

We start with some preparatory results. The next proposition shows that locally around its singular locus, any normal crossing variety is a Cartier divisor of a smooth variety.

\begin{proposition}[~\cite{Tzi09}]\label{-1}
Let $X$ be a scheme with normal crossing singularities. Then there is an analytic neighborhood $U$ of its singular locus and an embedding $U \subset V$, where $V$ is smooth and $\dim V = \dim U +1$.
\end{proposition}

The next proposition shows how to calculate $T^1(X)$ from an embedding $X \subset Y$ of $X$ as a Cartier divisor of a smooth variety $Y$ and it is the key in obtaining a formula for $T^1(X)$.

\begin{proposition}\label{0}
Let $X$ be a scheme with normal crossing singularities and let $X \subset Y$ be an embedding such that $Y$ is smooth and $\dim Y = \dim X+1$. Then
\[
T^1(X)=\mathcal{O}_Y(X)\otimes \mathcal{O}_Z
\]
where $Z$ is the singular locus of $X$.
\end{proposition}

\begin{proof}
Dualizing the exact sequence
\[
0 \la \mathcal{O}_Y(-X)\otimes \mathcal{O}_X \la \Omega_Y \otimes \mathcal{O}_X \la \Omega_X \la 0
\]
we get the exact sequence
\[
\mathcal{O}_Y(X)\otimes \mathcal{O}_X \la \mathcal{E}xt^1_X(\Omega_X,\mathcal{O}_X) \la \mathcal{E}xt^1_Y(\Omega_Y\otimes \mathcal{O}_X,\mathcal{O}_X)=0
\]
Moreover, $T^1(X)=\mathcal{E}xt^1_X(\Omega_X,\mathcal{O}_X) $ is locally free of rank 1 on $Z$. Hence there is a surjection $\mathcal{O}_Y(X)\otimes \mathcal{O}_Z \la T^1(X)$ and therefore $T^1(X) \cong \mathcal{O}_Y(X)\otimes \mathcal{O}_Z$, as claimed.
\end{proof}

\begin{definition}\label{definition1}
Let $X$ be a scheme with only double point normal crossing singularities. Let $Z \subset X$ be the singular locus of $X$. Then we denote by $L(Z,X)$ the invertible sheaf $\omega_X \otimes \omega_Z^{-1}$.
\end{definition}

\begin{lemma}\label{0.1}
Let $X$ be a scheme with only normal crossing double point singularities and let $Z$ be its singular locus. Then 
\begin{enumerate}
\item  $L(Z,X)$ is a 2-torsion invertible sheaf on $Z$, i.e., $L^{\otimes 2}\cong \mathcal{O}_Z$. Moreover, let $\pi \colon \tilde{X} \la X$ be the normalization of $X$, $\tilde{Z}=\pi^{-1}(Z)$ and $\pi \colon \tilde{Z} \la Z$ the induced map. Then $\pi$ is \'etale and $\pi^{\ast}L \cong \mathcal{O}_{\tilde{Z}}$.
\item Let $X\subset Y$ be an embedding such that $Y$ is smooth and $\dim Y=\dim X+1$. Then
\[
T^1(X)=L(Z,X) \otimes \bigwedge^2 \mathcal{N}_{Z/Y}
\]
\end{enumerate} 
\end{lemma}
\begin{remark}
If $X$ is simple normal crossing, then it is not hard to see that $\pi \colon \tilde{Z} \la Z$ is trivial and hence $\pi_{\ast} \mathcal{O}_{\tilde{Z}}=\mathcal{O}_Z \oplus \mathcal{O}_Z$. From this it easily follows that $L(Z,X)=\mathcal{O}_Z$. However, in general this is not true.
\end{remark}
\begin{proof}[Proof of Lemma~\ref{0.1}]
By Proposition~\ref{-1}, $T^1(X)=\mathcal{O}_Y(X)\otimes \mathcal{O}_Z$. Standard adjunctions give that 
\[
\mathcal{O}_Y(X) \otimes \mathcal{O}_Z = \omega_Y^{-1}\otimes \omega_X \otimes \mathcal{O}_Z=\bigwedge^2 \mathcal{N}_{Z/Y} \otimes (\omega_X \otimes \omega_Z^{-1})
\]
It remains to show that $L(Z,X)$ is a 2-torsion invertible sheaf on $Z$. Let $\pi \colon \tilde{X} \la X$ be the normalization of $X$ and $\tilde{Z}=\pi^{-1}(Z)$. Then a local calculation shows that $\pi \colon \tilde{Z} \la Z$ is \'etale. Moreover, $\pi^{\ast} \omega_X =\omega_{\tilde{X}} \otimes \mathcal{O}_{\tilde{X}}(\tilde{Z})$. Hence
\[
\omega_{\tilde{Z}}= \omega_{\tilde{X}} \otimes \mathcal{O}_{\tilde{X}}(\tilde{Z})\otimes \mathcal{O}_{\tilde{Z}}=\pi^{\ast}(\omega_X \otimes \mathcal{O}_Z)
\]
Moreover, since $\tilde{Z} \la Z$ is \'etale, $\pi^{\ast}\omega_Z = \omega_{\tilde{Z}}$. Therefore
\[
\pi^{\ast}\omega_Z=\pi^{\ast}(\omega_X \otimes \mathcal{O}_Z).
\]
Hence $\pi^{\ast}L(Z,X) \cong \mathcal{O}_{\tilde{Z}}$ and consequently $L(Z,X) \otimes \pi_{\ast} \mathcal{O}_{\tilde{Z}} \cong \pi_{\ast} \mathcal{O}_{\tilde{Z}}$. Therefore, since $\pi_{\ast}\mathcal{O}_{\tilde{Z}}$ is a rank 2 locally free sheaf on $Z$, it follows that $L(Z,X)^{\otimes 2} \cong \mathcal{O}_Z$ as claimed.
\end{proof}

\begin{theorem}\label{1}
Let $X$ be a scheme with only double point normal crossing singularities. Let $Z\subset X$ be its singular part, $\pi \colon \tilde{X} \la X$ its normalization and $\tilde{Z}=\pi^{-1}(Z)$. Then, $\tilde{X}$, $\tilde{Z}$ and $Z$ are smooth, $\pi \colon \tilde{Z} \la Z$ is \'etale, $T^1(X)$ is a line bundle on $Z$ and
\[
T^1(X)=L(Z,X) \otimes \bigwedge^2\pi_{\ast}\mathcal{O}_{\tilde{Z}}(\tilde{Z})
\]
where $L(Z,X)$ is a 2-torsion invertible sheaf on $Z$, as stated in Lemma~\ref{0.1}.
\end{theorem}
\begin{remark}
If $X=X_1 \cup X_2$ is simple normal crossing with two irreducible components, then the above formula is just the familiar $T^1(X)=\mathcal{N}_{Z/X_1}\otimes \mathcal{N}_{Z/X_2}$~\cite{PiPe83}.
\end{remark}
\begin{proof}[Proof of Theorem~\ref{1}]
The result is local around the singular locus of $X$ so we may assume by Proposition~\ref{-1} that $X \subset Y$, where $Y$ is smooth of dimension $\dim X+1$. The normalization $\tilde{X}$ of $X$ is simply the blow up of $X$ along $Z$. Let $f \colon \tilde{Y} \la Y$ be the blow up of $Y$ along $Z$. Then $\tilde{X}=f_{\ast}^{-1}X$. Let $\pi \colon \tilde{X} \la X$ be the restriction of $f$ on $\tilde{X}$. Let $E\subset Y$ be the $f$-exceptional divisor. Then $E=\mathbb{P}_Z(\mathcal{N}_{Z/X}^{-1})$, $\tilde{Z}=E\cdot \tilde{X} =\tilde{X} \cap E$. Hence
\[
\mathcal{O}_{\tilde{Z}}(\tilde{Z})=\mathcal{O}_{\tilde{X}}(\tilde{Z}) \otimes \mathcal{O}_{\tilde{Z}}=
\mathcal{O}_{\tilde{Y}}(E)\otimes \mathcal{O}_{\tilde{Z}}=\mathcal{O}_E(E)\otimes \mathcal{O}_{\tilde{Z}}= \mathcal{O}_E(-1)\otimes \mathcal{O}_{\tilde{Z}}
\]
Now $\mathcal{O}_E(-\tilde{Z})=\mathcal{O}_E(k)\otimes f^{\ast} \mathcal{L}$, for some $k \in \mathbb{Z}$ and $\mathcal{L} \in \mathrm{Pic}(Z)$. Since
$\tilde{Z} \la Z$ has degree $2$, it follows that $k=-2$. Then
\[
f^{\ast} \mathcal{L}=\mathcal{O}_E (-\tilde{Z}) \otimes \mathcal{O}_E (-2E)
\]
By Proposition~\ref{0}, $T^1(X)=\mathcal{O}_Y(X)\otimes\mathcal{O}_Z$. Also,
\[
f^{\ast}\mathcal{O}_Y(X)=\mathcal{O}_{\tilde{Y}}(\tilde{X})\otimes \mathcal{O}_{\tilde{Y}}(2E)
\]
and hence $f^{\ast}T^1(X)=\mathcal{O}_E (\tilde{Z}) \otimes \mathcal{O}_E (2E)$ which implies that $\mathcal{L}=T^1(X)^{-1}$. Therefore,
\[
\mathcal{O}_E(-\tilde{Z})=\mathcal{O}_E(-2)\otimes f^{\ast} T^1(X)^{-1}
\]
Now the exact sequence
\[
0 \la \mathcal{O}_E(-3)\otimes f^{\ast} T^1(X)^{-1} \la \mathcal{O}_E(-1) \la  \mathcal{O}_E(-1)\otimes \mathcal{O}_{\tilde{Z}}= \mathcal{O}_{\tilde{Z}}(\tilde{Z}) \la 0
\]
gives that
\begin{gather*}
\pi_{\ast}\mathcal{O}_{\tilde{Z}}(\tilde{Z}) = R^1f_{\ast}(\mathcal{O}_E(-3)\otimes f^{\ast}T^1(X)^{-1})= R^1f_{\ast}\mathcal{O}_E(-3) \otimes T^1(X)^{-1} =\\
(f_{\ast}\mathcal{O}_E(1))^{-1} \otimes \bigwedge^2\mathcal{N}_{Z/Y} \otimes T^1(X)^{-1}=\mathcal{N}_{Z/Y}\otimes \bigwedge^2\mathcal{N}_{Z/Y} \otimes T^1(X)^{-1}=\mathcal{N}_{Z/Y}\otimes L(Z,X)
\end{gather*}
since by Proposition~\ref{0.1}, $T^1(X)=\bigwedge^2\mathcal{N}_{Z/Y} \otimes L(Z,X)$. Hence
\begin{gather*}
T^1(X)=\bigwedge^2 \mathcal{N}_{Z/Y}\otimes L(Z,X)=(\bigwedge^2 (\pi_{\ast}\mathcal{O}_{\tilde{Z}}(\tilde{Z}) \otimes L(Z,X))) \otimes L(Z,X)= \\
(\bigwedge^2 (\pi_{\ast}\mathcal{O}_{\tilde{Z}}(\tilde{Z})) \otimes L(Z,X)^{\otimes 2} \otimes L(Z,X)=L(Z,X) \otimes \bigwedge^2  \pi_{\ast}\mathcal{O}_{\tilde{Z}}(\tilde{Z} )
\end{gather*}
as claimed.
\end{proof}

Next we extend Theorem~\ref{1} for any normal crossing variety. One of the differences of the general case is that the singular locus $Z$ of $X$ is no longer smooth. In fact it is not even Cohen-Macauley. For this reason it is preferable to work with a smooth model of $Z$. As a first step we shall rewrite the result of Theorem~\ref{1} in $\tilde{Z}$ instead of $Z$. In doing so we need the next result.

\begin{lemma}\label{1.1}
Let $f \colon X \la Y$ be a finite flat morphism of degree $d$ between smooth varieties. Let $\mathcal{E}$ be a locally free sheaf of rank $n$ on $E$. Then $f_{\ast}\mathcal{E}$ is locally free of rank $nd$ and
\[
f^{\ast}(\bigwedge^{nd}f_{\ast}\mathcal{E})=(\bigwedge^n \mathcal{E})^{\otimes d} \otimes f^{\ast} (\mathrm{det}(f_{\ast} \mathcal{O}_X))
\]
where $\mathrm{det}(f_{\ast} \mathcal{O}_X)=\wedge^d f_{\ast}\mathcal{O}_X$.
\end{lemma}

\begin{proof}
Let $A$ be an ample invertible sheaf on $Y$. Then $B=f^{\ast}A$ is ample on $X$. Then, since $X$ is smooth, $\mathcal{E}$ has a finite resolution
\[
0 \la \oplus_{i=1}^{k_s}B^{\nu_s} \la \cdots \la \oplus_{i=1}^{k_1}B^{ \nu_1} \la \mathcal{E} \la 0
\]
We will prove the lemma by doing induction on the minimum length $s(\mathcal{E})$ of such a resolution of $\mathcal{E}$. If $s(\mathcal{E})=1$, then $\mathcal{E}\cong \oplus_{i=1}^n B^{ k}$, for some $k \in \mathbb{Z}$. Then $f_{\ast}\mathcal{E} = \oplus_{i=1}^n (A^{k} \otimes f_{\ast}\mathcal{O}_X)$ and therefore
\[
\bigwedge^{nd}f_{\ast}\mathcal{E} = A^{ndk} \otimes \bigwedge^d f_{\ast}\mathcal{O}_X
\]
Hence
\[
f^{\ast} \bigwedge^{nd} f_{\ast}\mathcal{E} = B^{ndk}\otimes f^{\ast} (\mathrm{det}(f_{\ast} \mathcal{O}_X))=(\bigwedge^n \mathcal{E})^{\otimes d} \otimes f^{\ast} (\mathrm{det}(f_{\ast} \mathcal{O}_X))
\]
Now assume the result is true for all locally free sheaves $\mathcal{F}$ such that $s(\mathcal{F}) \leq k-1$. Let $\mathcal{E}$ be locally free with $s(\mathcal{E})=k$. Then there is an exact sequence
\[
0 \la \mathcal{F} \la \oplus_{i=1}^l B^m \la \mathcal{E} \la 0
\]
and therefore $\mathcal{F}$ is locally free with $s(\mathcal{F})=k-1$. Now after applying $f_{\ast}$, the induction hypotheses and some straightforward calculations we get the result.

\end{proof}

\begin{lemma}\label{1.2}
Let $X$ be a scheme with only double point normal crossing singularities and let $Z$ be its singular locus. Let $\pi \colon \tilde{X} \la X$ be the normalization and $\tilde{Z}=\pi^{-1}Z$. Then there is an exact sequence
\[
0 \la \mathcal{O}_Z \la \pi_{\ast} \mathcal{O}_{\tilde{Z}} \la L(Z,X) \la 0
\]
and $L(Z,X)$ is a 2-torsion invertible sheaf on $Z$, as stated in Lemma~\ref{0.1}.
\end{lemma}

\begin{proof}
We follow the main steps of the proof of Theorem~\ref{1}. Embed $X$ in a smooth variety $Y$ of dimension $\dim X +1$. Let $f \colon \tilde{Y} \la Y$ be the blow up of $Y$ along $Z$ and $E$ the $f$-exceptional divisor. Then $\tilde{X}=f_{\ast}^{-1}X$ and $\tilde{Z}=E\cdot \tilde{X}$. Then as we have seen in the proof of Theorem~\ref{1},
\[
\mathcal{O}_E(-\tilde{Z})=\mathcal{O}_E(-2)\otimes f^{\ast}(T^1(X)^{-1})
\]
Then the exact sequence
\[
0 \la \mathcal{O}_E(-\tilde{Z}) \la \mathcal{O}_E \la \mathcal{O}_{\tilde{Z}} \la 0
\]
give the exact sequence
\[
0 \la f_{\ast}\mathcal{O}_E=\mathcal{O}_Z \la f_{\ast} \mathcal{O}_{\tilde{Z}} \la R^1f_{\ast}(\mathcal{O}_E(-2)\otimes f^{\ast}(T^1(X)^{-1}) \la 0
\]
Moreover,
\begin{gather*}
R^1f_{\ast}(\mathcal{O}_E(-2)\otimes f^{\ast}(T^1(X)^{-1}) = R^1f_{\ast}\mathcal{O}_E(-2)\otimes T^1(X)^{-1} =\\
\bigwedge^2\mathcal{N}_{Z/Y}\otimes T^1(X)^{-1}=L(Z,X)
\end{gather*}
since by Lemma~\ref{0.1} $T^1(X)=L(Z,X) \otimes \bigwedge^2\mathcal{N}_{Z/Y}$, and $L(Z,X)=\omega_X \otimes \omega_Z^{-1}$ is a 2-torsion sheaf on $Z$.

\end{proof}

\begin{definition}
Let $X$ be a scheme with normal crossing singularities. We denote by $X^{max}$ the set of points having maximal multiplicity and by $X^{\geq s}$ the set of points of multiplicity at least $s$. A straightforward local calculation shows that $X^{max}$, $X^{\geq s}$ are closed subschemes of $X$ and that $X^{max}$ is smooth.
\end{definition}

\begin{theorem}\label{2}
Let $X$ be a scheme with normal crossing singularities and let $Z$ be its singular locus. Construct inductively the sequence of morphisms
\[
X^{\prime}=X_k \stackrel{f_k}{\la} X_{k-1} \stackrel{f_{k-1}}{\la} \cdots \stackrel{f_2}{\la} X_1 \stackrel{f_1}{\la} X_0=X
\]
such that $X_{i+1} $ is the blow up of $X_i$ along $X_i^{max}$. Let $f=f_1 \circ \cdots \circ f_k$. Let $Z^{\prime} \subset X^{\prime}$ be the divisorial part of $f^{-1}(Z)$ that dominates $Z$ and $E_s$ the reduced $f$-exceptional divisor that dominates $X^{\geq s}$, $s\geq 3$. Then $X^{\prime}$ and $Z^{\prime}$ are smooth and
\[
f^{\ast}T^1(X)=\mathcal{O}_{Z^{\prime}}(Z^{\prime})^{\otimes 2} \otimes (\otimes_{s \geq 3} \mathcal{O}_{X^{\prime}}(sE_s))
\]
\end{theorem}

\begin{proof}
Embed $X$ in a smooth $n+1$-dimensional variety $Y$, where $n=\dim X$. We proceed in two steps.

\textbf{Step 1.} First we succesively blow up the locus of points of multiplicity $\geq 3$ in order to reduce the calculation to that of a normal crossing scheme with only double points. Let $m$ be the maximal multiplicity of the singularities of $X$. If $m=2$, then go to step 2. Suppose that $m\geq 3$. Then locally at a point of maximal multiplicity, $X$ is given by $x_1\cdots x_m=0 \subset \mathbb{C}^{n+1}$, and $X^{max}$ by $x_1=x_2=\cdots =x_m=0$. Let $f_1 \colon Y_1 \la Y$ be the blow up of $Y$ along $X^{max}$ and $F_1$ the $f_1$-exceptional divisor. Let $X_1 =(f_1)_{\ast}^{-1}X$, $Z_1=(f_1)_{\ast}^{-1}Z$ and $B_1=F_1 \cdot X_1$. A straightforward local calculation shows that $X_1$ has normal crossing singularities of maximal multiplicity $m_1=m-1$, $B_1$ is $f_1$-exceptional, it is not contained in the singular locus of $X_1$ and that $X_1^{max} \not \subset B_1$. Moreover,
\begin{equation}\label{eq2.1}
f_1^{\ast} \mathcal{O}_Y(X)=\mathcal{O}_{Y_1}(X_1) \otimes \mathcal{O}_{Y_1}(mB_1)
\end{equation}
Also, by Proposition~\ref{0}, $T^1(X)=\mathcal{O}_Y(X)\otimes \mathcal{O}_Z$. Therefore~(\ref{eq2.1}) gives that,
\[
f_1^{\ast}T^1(X)=T^1(X_1)\otimes \mathcal{O}_{Y_1}(mB_1)
\]
Repeating the previous process of blowing up the locus of highest multiplicity, we get a sequence of maps
\[
X_{m-2}\stackrel{f_{m-2}}{\la} \cdots \stackrel{f_{2}}{\la} X_1 \stackrel{f_{1}}{\la} X
\]
where $X_{m-2}$ has only normal crossing points of multiplicity 2. By a slight abuse of notation, denote by $B_i$ the birational transform of the $f_i$ exceptional divisor in $X_{m-2}$ and let $g = f_{m-2}\circ \cdots \circ f_1$. Then,
\begin{equation}\label{eq2.2}
g^{\ast}T^1(X)=T^1(X_{m-2})\otimes (\otimes_{i=1}^{m-2} \mathcal{O}_{X^{\prime}}((m-i+1)B_i))
\end{equation}

\textbf{Step 2.} Let $f_{m-1} \colon X_{m-1} \la X_{m-2}$ be the normalization of $X_{m-2}$, which is smooth since $X_{m-2}$ has only double point singularities. Let $Z_{m-2}$ its singular locus and $Z_{m-1}=f_{m-1}^{-1}(Z_{m-2})$. Then by Theorem~\ref{1},
\[
T^1(X_{m-2})= L(Z_{m-2},X_{m-2}) \otimes \bigwedge^2 (f_{m-1})_{\ast}\mathcal{O}_{Z_{m-1}}(Z_{m-1})
\]
and by Lemmas~\ref{0.1},~\ref{1.1},~\ref{1.2},
\begin{equation}\label{eq2.3}
f_{m-1}^{\ast}T^1(X_{m-2})=\mathcal{O}_{Z_{m-1}}(Z_{m-1})^{\otimes 2}
\end{equation}
Let $f=f_{m-1}\circ g$ and $B^{\prime}_i=f_{m-1}^{-1}B_i$. Then from (\ref{eq2.2}),~(\ref{eq2.3}) it follows that
\begin{equation}\label{eq2.4}
f^{\ast}T^1(X)=\mathcal{O}_{Z_{m-1}}(Z_{m-1})^{\otimes 2} \otimes (\otimes_{i=1}^{m-2} \mathcal{O}_{X^{\prime}}((m-i+1)B^{\prime}_i))
\end{equation}
Note that by construction, $B^{\prime}_i$ dominates the locus of points of multiplicity $\geq m-i+1$. Now setting $X^{\prime}=X_{m-1}$, $E_{m-i+1}=B_i $,~(\ref{eq2.4}) takes the form stated in the theorem.
\end{proof}

\begin{remark}
One may try to get a formula for $T^1(X)$ in the normalization $\pi \colon \tilde{X} \la X$. However, $\tilde{Z}=\pi^{-1}(Z)$ is singular and it is preferable to work with smooth varieties. The pair $(Z^{\prime},X^{\prime})$ is a log resolution for $(Z,X)$ that is obtained in a natural way by repeatedly blowing up the centers of maximal multiplicity. The disadvantage of this approach is that $(Z^{\prime},X^{\prime})$ is not characterized by any numerical property that would make it unique, as for example the minimal log resolution in the case of surfaces. However, in dimension at most 3 we can get a more natural description by using the minimal model program.
\end{remark}

\begin{theorem}\label{3}
Let $X$ be a scheme with normal crossing singularities with $\dim X \leq 3$ and let $Z$ be its singular locus. Let $\phi \colon \tilde{X} \la X$ be a birational morphism such that
\begin{enumerate}
\item $(\tilde{X},1/2\tilde{Z})$, $\tilde{Z}$ are terminal;
\item $K_{\tilde{X}}+1/2\tilde{Z}$ and $K_{\tilde{Z}}$ are $\phi$-nef.
\end{enumerate}
where $\tilde{Z} \subset \tilde{X}$ be the reduced divisorial part of $\phi^{-1}(Z)$ that dominates $Z$ (such spaces do exist). Then
\[
\phi^{\ast}T^1(X)=\mathcal{O}_{\tilde{Z}}(\tilde{Z})^{\otimes 2} \otimes \mathcal{O}_{\tilde{X}}(3E)\otimes \mathcal{O}_{\tilde{Z}}
\]
where $E\subset \tilde{X}$ is the reduced $\phi$-exceptional divisor that dominates the set of singular points of $X$ of multiplicity at least three.
\end{theorem}
\begin{remark}
\begin{enumerate}
\item In the case of surfaces, $(\tilde{Z},\tilde{X})$ is simply the minimal log resolution of $(Z,X)$.
\item In the case of 3-folds,~(\ref{3}.1) implies that $\tilde{Z}$ is smooth in $\tilde{X}$. However, the form written seems to be more natural to generalize in higher dimensions.
\item The proof of the theorem shows that if $\overline{Z}$ is the normalization of $Z$, then the induced map $\tilde{Z} \la \overline{Z}$ is \'etale.
\item The problems in order to obtain a similar result in all dimensions are the following. The existence of a pair $(\tilde{X},\tilde{Z})$ with the properties stated in Theorem~\ref{3} is not a formal consequence of the minimal model program. One could start with a log resolution $(X^{\prime},Z^{\prime})$ of $(X,Z)$ and then try to run a MMP simultaneously for $(X^{\prime},1/2Z^{\prime})$ and $Z^{\prime}$, but this cannot be done in general. The construction of $(\tilde{X},\tilde{Z})$ is explicit in the proof of Theorem~\ref{3} and in principle if one is careful enough it should be possible to generalize the argument in all dimensions.
\end{enumerate}
\end{remark}

\begin{proof}[Proof of Theorem~\ref{3}]
We only do the 3-fold case. The surface case is much simpler.

The proof consists of two steps. In the first one we will explicitly construct a pair $(\tilde{X},\tilde{Z})$ with the properties of the statement and in the second part we will show that given any other pair $(X^{\prime},Z^{\prime})$ having the same properties, $T^1(X)$ is given by the same formula.

\textbf{Step 1.} Let $X \subset Y$ be an embedding of $X$ into a smooth 4-fold. Since $\dim X=3$, $\mathrm{mult}_P(X) \leq 4$, for all $P\in X$. Now repeat the construction in Theorem~\ref{2}. There is a sequence of maps
\[
X^{\prime} \stackrel{g}{\la} X_2 \stackrel{f_2}{\la} X_1 \stackrel{f_1}{\la} X
\]
with the following properties. $X_1$ is the blow up of $X$ along the locus of points of multiplicity $4$. Then the singular points of $X_1$ have multiplicity at most $3$ and $X_2$ is the blow up of $X_1$ along the locus of points of multiplicity $3$. $X_2$ has only double points and $X^{\prime}$ is its normalization.

Let $E_1$ be the $f_1$-exceptional divisor and $Z_1=(f_1)_{\ast}^{-1}Z$. A straightforward local calculation shows that $Z_1$ is the singular locus of $X_1$,  and over any singular point of $X$ with multiplicity $4$, $E_1 \cong (x_0x_1x_2x_3=0)\subset \mathbb{P}^3$. Moreover, $X_1^{max} \cap E_1 = \{P_1,P_2,P_3,P_4\}$, where $P_i$ are the vertices of the tetrahedron $x_0x_1x_2x_3=0$ in $\mathbb{P}^3$, and $E_1 \cap Z_1 =\cup_{i,j} L_{i,j}$, where $L_{i,j}$ is the line connecting the verices $P_i$ and $P_j$.

Let $E_2$ be the $f_2$-exceptional divisor and $Z_2=(f_2)_{\ast}^{-1}(Z_1)$. Let $P\in X_1$ be a point such that $\mathrm{mult}_P(X_1)=3$. Then a straightforward local calculation shows that $X_2$ has only double points, $f_2^{-1}(P)\cong(x_0x_1x_2=0)\subset \mathbb{P}^2$, $Z_2$ is the singular locus of $X_2$ and $Z_2 \cap f_2^{-1}(P) =\{Q_1,Q_2,Q_3\}$, where $Q_i$ are the vertices of the triangle $x_0x_1x_2=0$ in $\mathbb{P}^3$.

Let $\tilde{E}_1=(f_2)_{\ast}^{-1}E_1$. Then $\tilde{E}_1=f_2^{\ast}E_1$. Moreover, $\tilde{E}_1$ is the blow up of $E_1$ along its vertices, and over any point of multiplicity $4$ of $X$, $\tilde{E}_1 \cap E_2 = \cup_{i=1}^4 f_2^{-1}(P_i)$. Moreover, from the previous discussion it follows that $f_2^{-1}(P_i)=F_{i,1}\cup F_{i,2}\cup F_{i,3}$, where $F_{i,j}$ are the edges of the triangle $(x_0x_1x_2=0) \subset \mathbb{P}^2$, and $\tilde{E}_1 \cap Z_2 = \cup_{i,j} \tilde{L}_{i,j}$, where $\tilde{L}_{i,j}$ are the birational transforms of the edges $L_{i,j}$ of $E_1$ in $E_2$.

The normalization $X^{\prime}$ of $X_2$ is the blow up of $X_2$ along $Z_2$. Let $Z^{\prime}=g^{-1}(Z_2)$. Then $g \colon Z^{\prime} \la Z_2$ is \'etale. Denote by $D^{\prime}$ the birational transform of any divisor $D \subset X_2$ in $Z^{\prime}$. Then $E_2^{\prime}$, $E_1^{\prime}$ are the blow ups of $E_2$, $E_1$ along their edges and a straightforward local calculation shows that $E_2^{\prime}$, $E_1^{\prime}$ are smooth. Moreover, over any singular point $P\in X$ of multiplicity $4$, $E_2^{\prime}=\coprod_{i=1}^4 H_i $, where $H_i$ is isomorphic to the blow up of $\mathbb{P}^2$ along the vertices of the triangle $x_0x_1x_2=0$. In other words, $E^{\prime}_1$ is the log resolution of $E_1$ and its edges obtained by first blowing up the vertices and then the edges. Hence in each $H_i$ there is the following configuration of $f$-exceptional curves, where $f=f_1 \circ f_2 \circ g$.
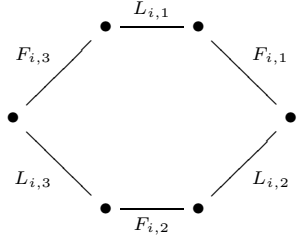
\[
\xymatrix{
              & \bullet   \ar@{-}[dl]_{F_{i,3}}\ar@{-}[r]^{L_{i,1}}                  &        \bullet \ar@{-}[dr]^{F_{i,1}}             &               \\
\bullet \ar@{-}[dr]_{L_{i,3}}       &                            &                             &      \bullet \ar@{-}[dl]^{L_{i,2}}            \\
              &     \bullet \ar@{-}[r]_{F_{i,2}}                &         \bullet              &        \\
}
\]
Moreover, we have the following intersection table.
\begin{gather}\label{intersection-table}
L_{i,j}^2=F_{i,j}^2=-1 \notag \\
\cup_{i,j}L_{i,j} = H_i \cdot Z^{\prime}  \\
\cup_{i,j} F_{i,j} = H_i \cdot E_2^{\prime} \notag \\
L_{i,j} \cdot E^{\prime}_2=F_{i,j}\cdot Z^{\prime} =2 \notag \\
L_{i,j}\cdot E_1^{\prime}=L_{i,j} \cdot (f_2\circ g)^{\ast} E_1 = (f_2\circ g)_{\ast}(L_{i,j}) \cdot E_1 =-1 \notag 
\end{gather}

\begin{gather*}
F_{i,j} \cdot E_1^{\prime} =F_{i,j} \cdot (f_2\circ g)^{\ast} E_1 = (f_2\circ g)_{\ast}(F_{i,j}) \cdot E_1=0  \\
F_{i,j} \cdot E^{\prime}_2 =-1  \\
L_{i,j} \cdot Z^{\prime} = (L_{i,j})^2_{E_1^{\prime}} =-1  \\
(L_{i,j})^2_{Z^{\prime}}=L_{i,j} \cdot E_1^{\prime}=-1 
\end{gather*}
where by $(L_{i,j})^2_{E_1^{\prime}}$, $(L_{i,j})^1_{Z^{\prime}}$ we denote the self-intersection numbers of $L_{i,j}$ in $E_1^{\prime}$ and $Z^{\prime}$, respectively.

Next we will use the Minimal Model Program to obtain a 3-fold with the properties stated in the theorem.

Standard adjunction formulas (which could be calculated using an embedding $X \subset Y$ of $X$ in a smooth 4-fold $Y$ and following the argument of the proof of Theorem~\ref{1}), give that
\[
K_{X^{\prime}}+1/2Z^{\prime} =f^{\ast}K_X-E_1^{\prime}-E_2^{\prime}-1/2Z^{\prime}.
\]
and therefore
\begin{gather*}
(K_{X^{\prime}}+1/2Z^{\prime})\cdot L_{i,j}=-1/2 \\
K_{X^{\prime}}\cdot L_{i,j} =0
\end{gather*}
Moreover, the exact sequence
\[
0 \la \mathcal{N}_{L_{i,j}/E_1^{\prime}} \la \mathcal{N}_{L_{i,j}/X^{\prime}} \la \mathcal{O}_{X^{\prime}}(E^{\prime}_1) \otimes \mathcal{O}_{L_{i,j}} \la 0
\]
give that
\[
\mathcal{N}_{L_{i,j}/X^{\prime}} = \mathcal{O}_{\mathbb{P}^1}(-1) \oplus \mathcal{O}_{\mathbb{P}^1}(-1)
\]
Hence there is a flopping contraction $h^{\prime} \colon X^{\prime} \la S$ contracting each $L_{i,j}$ to an ordinary 3-fold double point. Let $\psi \colon  X^{\prime} \dasharrow X^{\prime\prime} $ be its flop. This is a standard flop and its construction is described in the following diagram:
\begin{equation}\label{flop}
\xymatrix{
              &             W \ar[dl]_{\psi^{\prime}}\ar[dr]^{\psi^{\prime\prime}}               &              \\
X^{\prime} \ar[dr]_{h^{\prime}} \ar@{-->}[rr]^{\psi}         &                            &     X^{\prime\prime} \ar[dl]^{h^{\prime\prime}}\\
                                         &              S                 &      \\
}
\end{equation}
Here $W$ is the blow up of $X^{\prime}$ along $L_{i,j}$. Let $B$ be the $\psi^{\prime}$-exceptional divisor. Then over a neighborhood of any of the $L_{i,j}$, $B$ is a ruled surface over $L_{i,j}$ and in fact $B \cong \mathbb{P}^1 \times \mathbb{P}^1$. Then $\psi^{\prime\prime}$ is the contraction of the other ruling and hence $X^{\prime\prime}$ is also smooth. Let $L_{i,j}^{\prime\prime}$ be the $h^{\prime\prime}$-exceptional curves and $Z^{\prime\prime}$ be the birational transform of $Z^{\prime}$ in $X^{\prime\prime}$. Then from diagram~(\ref{flop}) and the explicit description of $X^{\prime}$, $E_1^{\prime}$ and $Z^{\prime}$ it follows that $Z^{\prime}_W=(\psi^{\prime})_{\ast}^{-1} Z^{\prime} \cong Z^{\prime}$ and therefore
\[
\phi \colon Z^{\prime} \la Z^{\prime\prime}
\]
is a morphism and is in fact, from~(\ref{intersection-table}), it is the contraction of the $(-1)$-curves of $Z^{\prime}$. Then $L^{\prime\prime}_{i,j} \cdot Z^{\prime\prime}=1$, $K_{X^{\prime\prime}}\cdot L^{\prime\prime}_{i,j}=0$ and therefore
\[
(K_{X^{\prime\prime}}+1/2 Z^{\prime\prime})=1/2
\]
However, $K_{X^{\prime\prime}}+1/2 Z^{\prime\prime}$ is still not nef. Let $F^{\prime\prime}_{i,j}$ be the birational transform of $F_{i,j}$ in $X^{\prime\prime}$. Then, since
\[
(\psi^{\prime\prime})^{\ast}Z^{\prime\prime}=(\psi^{\prime})_{\ast}^{-1}Z^{\prime}=(\psi^{\prime})^{\ast}Z^{\prime}-B
\]
it follows from~(\ref{intersection-table}) that
\[
Z^{\prime\prime}\cdot F^{\prime\prime}_{i,j} = Z^{\prime}\cdot F_{i,j} - B \cdot (\psi^{\prime})_{\ast}^{-1} F_{i,j} =2-2=0
\]
Moreover, $K_{X^{\prime\prime}} \cdot F^{\prime\prime}_{i,j} =K_{X^{\prime}}\cdot F_{i,j} = -1$ and hence
\[
(K_{X^{\prime\prime}}+1/2 Z^{\prime\prime}) \cdot F^{\prime\prime}_{i,j} = -1
\]
Let $H^{\prime\prime}_i=\phi_{\ast}^{-1} H_i$. Then from diagram~\ref{flop} it follows that $H_i^{\prime\prime} \cong \mathbb{P}^2$. Let us elaborate more on this. As mentioned earlier in the proof, if $P \in X$ is any point of multiplicity $4$ of $X$, then $f_1^{-1}(P)$ is isomorphic to the tetrahedron $x_0x_1x_2x_3=0$ in $\mathbb{P}^3$. Let $\Delta_i$ be its faces, $i=1,2,3,4$. Then $\Delta_i \cong \mathbb{P}^2$ and the edges $L_{i,j}$ of $f_1^{-1}(P)$ on $\Delta_i$ is the triangle $x_0x_1x_2=0$. Then $H_i$ is the blow up of $\Delta_i$ along the vertices of the triangle and $F_{i,j}$ are the exceptional curves. Now $H^{\prime\prime}_{i}$ is the contraction of the $L_{i,j}$. But this is exactly the construction of the standard quadratic transformation of $\mathbb{P}^2$. hence $H_i^{\prime\prime} \cong \mathbb{P}^2$ and
\[
\Delta_i=\mathbb{P}^2 \dasharrow \mathbb{P}^2\cong H_i^{\prime\prime}
\]
is the standard quadratic transformation of $\mathbb{P}^2$.

Hence $E_1^{\prime\prime}$ is a disjoint union of projective planes. Moreover,
\[
(\psi^{\prime\prime})^{\ast}E_1^{\prime\prime}=(\psi^{\prime})^{\ast}E_1^{\prime}-B
\]
and hence
\[
F_{i,j}^{\prime\prime} \cdot E_1^{\prime\prime}=F_{i,j}\cdot E_1^{\prime}-F_{i,j}\cdot B =0-2=-2
\]
Therefore, there is a birational morphism $\alpha \colon X^{\prime\prime} \la \tilde{X} $ over $X$, contractiong every irreducible component of $E_1^{\prime\prime}$ to a cyclic quotient singularity of type $\frac{1}{2}(1,1,1)$. let $\tilde{Z}=\alpha_{\ast}Z^{\prime\prime}$. Then $\tilde{Z} \cong Z^{\prime\prime}$ and from the construction it is also clear that it is contained in the smooth part of $\tilde{X}$. Moreover, $K_{\tilde{X}}+1/2\tilde{Z}$ and $K_{\tilde{Z}}$ are both nef over $X$. Hence $(\tilde{X},\tilde{Z})$ satisfies the numerical properties of the theorem. It is also clear from the above construction that the natural induced map $\tilde{Z} \la \overline{Z}$, where $\overline{Z}$ is the normalization of $Z$, is \'etale.

Next we show that $T^1(X)$ is given by the formula claimed in the statement. We do it by moving around the diagram
\begin{equation}\label{construction-diagram1}
\xymatrix{
              &             W \ar[dl]_{\psi^{\prime}}\ar[dr]^{\psi^{\prime\prime}}               &              \\
X^{\prime} \ar[d]_{f} \ar@{-->}[rr]^{\psi}         &                            &     X^{\prime\prime} \ar[d]^{\alpha}\\
                  X                       &                                                  &   \tilde{X} \ar[ll]^{\phi}    \\
}
\end{equation}
and the corresponding for $Z$,
\begin{equation}\label{construction-diagram2}
\xymatrix{
Z^{\prime} \ar[rr]^{\psi}\ar[dr]_{f} & & Z^{\prime\prime}=\tilde{Z} \ar[dl]^{\phi} \\
                                     & Z &
}
\end{equation}
Here we should remark that from the construction of diagram~(\ref{construction-diagram1}) it follows that $(\psi^{\prime})_{\ast}^{-1}Z^{\prime} \cong Z^{\prime}$ and hence $\psi \colon Z^{\prime} \la Z^{\prime\prime}$ is just $\psi^{\prime\prime} \colon (\psi^{\prime})_{\ast}^{-1}Z^{\prime} \la Z^{\prime\prime}$.

From Theorem~\ref{2}, it follows that
\begin{equation}\label{eq-3.1}
f^{\ast}T^1(X)=\mathcal{O}_{Z^{\prime}}(Z^{\prime})^{\otimes 2} \otimes \mathcal{O}_{X^{\prime}}(4E_1^{\prime}) \otimes \mathcal{O}_{X^{\prime}}(3E_2^{\prime})
\otimes \mathcal{O}_{Z^{\prime}}
\end{equation}
Next we reduce it to $Z^{\prime\prime}=\tilde{Z}$. First we get a formula in $X^{\prime\prime}$. From diagram~(\ref{construction-diagram1}) and standard adjunctions we get that
\begin{gather*}
\psi^{\ast}\mathcal{O}_{Z^{\prime\prime}}(2Z^{\prime\prime})=\mathcal{O}_{Z^{\prime}}(2Z^{\prime})\otimes \mathcal{O}_W(-2B) \otimes \mathcal{O}_{Z^{\prime}}\\
\psi^{\ast} ( \mathcal{O}_{X^{\prime\prime}}(E_2^{\prime\prime})\otimes \mathcal{O}_{Z^{\prime\prime}})=( \mathcal{O}_{X^{\prime}}(E_2^{\prime})\otimes \mathcal{O}_{Z^{\prime}})\otimes \mathcal{O}_W(2B) \otimes \mathcal{O}_{Z^{\prime}}\\
\psi^{\ast} ( \mathcal{O}_{X^{\prime\prime}}(E_1^{\prime\prime})\otimes \mathcal{O}_{Z^{\prime\prime}})=( \mathcal{O}_{X^{\prime}}(E_1^{\prime})\otimes \mathcal{O}_{Z^{\prime}})\otimes \mathcal{O}_W(-B) \otimes \mathcal{O}_{Z^{\prime}}\\
\end{gather*}
Hence
\[
(\phi \circ \alpha)^{\ast}T^1(X)=\mathcal{O}_{Z^{\prime\prime}}(Z^{\prime\prime})^{\otimes 2} \otimes \mathcal{O}_{X^{\prime\prime}}(4E_1^{\prime\prime}) \otimes \mathcal{O}_{X^{\prime\prime}}(3E_2^{\prime\prime})
\otimes \mathcal{O}_{Z^{\prime}}
\]
and so the formula is unchanged in $X^{\prime\prime}$ (we did not expect a change since $X^{\prime}$ and $X^{\prime\prime}$ are isomorphic in codimension 1). Moreover, a carefull look at the construction of $X^{\prime\prime}$ reveals that $E_1^{\prime\prime} \cap Z^{\prime\prime} = \emptyset$. Hence
\[
(\phi \circ \alpha)^{\ast}T^1(X)=\mathcal{O}_{Z^{\prime\prime}}(Z^{\prime\prime})^{\otimes 2} \otimes \mathcal{O}_{X^{\prime\prime}}(3E_2^{\prime\prime})
\otimes \mathcal{O}_{Z^{\prime}}
\]
and therefore
\[
\phi^{\ast}T^1(X)=\mathcal{O}_{\tilde{Z}}(\tilde{Z})^{\otimes 2} \otimes \mathcal{O}_{\tilde{X}}(3E)
\otimes \mathcal{O}_{\tilde{Z}}
\]
where $E=\alpha_{\ast}E_2^{\prime\prime}$, as claimed in the statement.

\textbf{Step 2.} Let $\psi \colon (\hat{X},\hat{Z}) \la (X,Z)$ be a morphism such that $(\hat{X},1/2\hat{Z})$ is terminal and  $K_{\hat{X}}+1/2\hat{Z}$, $K_{\hat{Z}}$ are $\psi$-nef. We will show that $\psi^{\ast} T^1(X)$ is given by the formula stated in the theorem.

Since $(\hat{X},1/2\hat{Z})$ is terminal and nef over $X$, it follows that there is a birational map
\[
g \colon \tilde{X} \dasharrow \hat{X}
\]
which is an isomorphism in codimension 1. Moreover, since $K_{\hat{Z}}$ and $K_{\tilde{Z}}$ are nef over $Z$, $g$ induces an isomorphism between $\tilde{Z}$ and $\hat{Z}$. Since we already know that $\tilde{Z} \la \overline{Z}$ is \'etale, $\tilde{Z}$ does not contain any $\phi$-exceptional curves. Therefore $\hat{Z}$ does not contain any $\psi$-exceptional curves and hence there is a commutative diagram
\[
\xymatrix{
\mathrm{Ref}(\tilde{X}) \ar[r]^{g_{\ast}}\ar[d] & \mathrm{Ref}(\hat{X}) \ar[d] \\
\mathrm{Pic}(\tilde{Z}) \ar[r]^{g_{\ast}} & \mathrm{Pic}(\hat{Z})
}
\]
where $\mathrm{Ref}(\tilde{X})$, $\mathrm{Ref}(\hat{X})$ are the groups of rank 1 reflexive sheaves on $\tilde{X}$ and $\hat{X}$, respectively. Moreover, the horizontal maps are isomorphisms and the vertical the restrictions (the restriction maps are well defined because $\tilde{X}$, $\hat{X}$ are terminal and $\tilde{Z}$, $\hat{Z}$ smooth).

Since $g$ is an isomorphism in codimension 1, $\phi$ and $\psi$ have the same exceptional divisors. Hence if $\tilde{E}$ is $\phi$-exceptional dominating the locus of points of multiplicity at least 3, $\hat{E}=g_{\ast}^{-1} \tilde{E}$ is $\psi$-exceptional dominating the set of points of multiplicity at least 3 and $g_{\ast}\mathcal{O}_{\tilde{X}}(\tilde{E})=\mathcal{O}_{\hat{X}}(\hat{E})$. Moreover, $g_{\ast}\mathcal{O}_{\tilde{X}}(\tilde{Z})=\mathcal{O}_{\hat{X}}(\hat{Z})$. Hence
\[
\psi^{\ast}T^1(X)=g_{\ast}\phi^{\ast}T^1(X)=g_{\ast}(\mathcal{O}_{\tilde{Z}}(\tilde{Z})^{\otimes 2} \otimes \mathcal{O}_{\tilde{X}}(3\tilde{E})
\otimes \mathcal{O}_{\tilde{Z}})=\mathcal{O}_{\hat{Z}}(\hat{Z})^{\otimes 2} \otimes \mathcal{O}_{\hat{X}}(3\hat{E})
\otimes \mathcal{O}_{\hat{Z}}
\]
as claimed.
\end{proof}

\section{Obstructions}

Let $X$ be a variety with normal crossing singularities. It is well known that $H^2(T_X)$ and $H^1(T^1(X))$ are obstruction spaces to lift deformations of $X$. In this section we describe these spaces and in fact we show that $H^2(T_X)=0$. First we present some preliminary results.

\begin{proposition}[~\cite{Fr83}]\label{ob1}
Let $X$ be a scheme with only normal crossing singularities. Let $\tau_X \subset \Omega_X$ be the torsion subsheaf of $\Omega_X$. Then
\begin{gather*}
\Omega_X/\tau_X \cong \Omega_X^{\ast\ast} \\
\mathrm{Ext}_X^i(\Omega_X/\tau_X,\mathcal{O}_X)\cong H^i(T_X) \\
\mathrm{Ext}_X^i(\tau_X,\mathcal{O}_X)=H^{i-1}(T^1(X))
\end{gather*}
\end{proposition}

\begin{corollary}\label{ob2}
Let $X$ be a projective scheme with only normal crossing singularities. Then
\[
H^2(T_X)=H^{n-2}((\Omega_X/\tau_X) \otimes \omega_X)
\]
where $\dim X =n$.
\end{corollary}

\begin{proof}
By Proposition~\ref{ob1}, $H^2(T_X)=\mathrm{Ext}_X^2(\Omega_X/\tau_X,\mathcal{O}_X)=H^{n-2}((\Omega_X/\tau_X) \otimes \omega_X)$, by Serre duality.
\end{proof}

The next result is the main technical tool in order to calculate $H^2(T_X)$.

\begin{proposition}\label{ob3}
Let $X$ be a variety with normal crossing singularities. Let $Z \subset X$ be its singular locus and $\pi \colon \tilde{X} \la X$, $\nu \colon \tilde{Z} \la Z$ be the normalization of $X$ and $Z$, respectively. Then there is an exact sequence
\begin{equation}\label{ob-seq}
0 \la \tau_X \la \Omega_X \stackrel{\delta_1}{\la} \pi_{\ast}\Omega_{\tilde{X}} \stackrel{\delta_2}{\la} \nu_{\ast}(\Omega_{\tilde{Z}} \otimes L)
\end{equation}
where $L$ is an invertible sheaf on $\tilde{Z}$ such that $L^{\otimes 2} \cong \mathcal{O}_{\tilde{Z}}$. In particular, if $X$ has at most double points, then $\delta_2$ is surjective and the sequence is exact on the right too.
\end{proposition}

\begin{remark}
If $X=\cup_{i=1}^NX_i$ is a reducible simple normal crossing variety, then $L \cong \mathcal{O}_Z$ and the exact sequence~(\ref{ob-seq}) is part of a long sequence~\cite{Fr83}. In the general case I believe there must be a long exact sequence
\[
0 \la \tau_X \la \Omega_X \la \pi_{\ast}\Omega_{\tilde{X}} \stackrel{\delta_1}{\la} \nu_{\ast}(\Omega_{\tilde{Z}}\otimes L) \la (\nu_1)_{\ast}(\Omega_{\tilde{Z_1}}\otimes L_1) \la (\nu_2)_{\ast}(\Omega_{\tilde{Z_2}}\otimes L_2) \la \cdots
\]
where the $Z_i$ are defined inductively by $Z_0=Z$, $Z_i$ is the singular locus of $Z_{i-1}$, $\nu_i \colon \tilde{Z}_i \la Z_i$ is the normalization and $L_i$ are 2-torsion invertible sheaves on $\tilde{Z}_i$. However, at the moment I am having a few technical difficulties to prove this. The current statement of proposition~\ref{ob3} is sufficient to treat the cases when $\dim X \leq 3$.
\end{remark}

\begin{proof}[Proof of Proposition~\ref{ob3}]

The proof consists of two steps. In the first one we consider the case when $X$ has only double points and in the second we do the general case.

\textbf{Step 1.} Suppose that $X$ has only double points. Then $\tilde{Z}=Z$. The idea of the proof is to first construct the maps $\delta_i$ in a natural unique way and then check the excatness locally.

Let $\hat{Z}=\pi^{-1}(Z)$. Then $\hat{Z}$ is smooth and $\hat{Z} \la Z$ is \'etale. Moreover, by Lemma~(\ref{1.2}), there is an exact sequence
\begin{equation}\label{eq-ob3-0}
0 \la \mathcal{O}_Z \stackrel{\beta}{\la} \pi_{\ast}\mathcal{O}_{\hat{Z}} \stackrel{\gamma}{\la} L \la 0
\end{equation}
where $L=L(Z,X)$ is a 2-torsion invertible sheaf on $Z$. The natural map $\pi^{\ast}\Omega_X \la \Omega_{\tilde{X}}$ gives a map $\pi_{\ast}\pi^{\ast}\Omega_X \la \pi_{\ast}\Omega_{\tilde{X}}$ and hence a map $\delta_0 \colon \Omega_X \la \pi_{\ast}\Omega_{\tilde{X}}$. Since $\tilde{X}$ is smooth, $\pi_{\ast}\Omega_{\tilde{X}}$ is torsion free and therefore the kernel of $\delta_0$ is the torsion of $\Omega_X$. Therefore there is an exact sequence
\[
0 \la \tau_X \la \Omega_X \stackrel{\delta_0}{\la} \pi_{\ast}\Omega_{\tilde{X}}
\]
Next we define $\delta_1 \colon \pi_{\ast}\Omega_{\tilde{X}} \la \nu_{\ast}(\Omega_{\tilde{Z}}\otimes L)=\Omega_Z\otimes L$. There is a natural surjection $\Omega_{\tilde{X}}\la \Omega_{\hat{Z}}$ and therefore a surjection
\begin{equation}\label{eq-ob3-1}
\pi_{\ast}\Omega_{\tilde{X}} \la \pi_{\ast} \Omega_{\hat{Z}} \la 0
\end{equation}
Now $\hat{Z} \la Z$ is \'etale and hence $\pi^{\ast}\Omega_Z = \Omega_{\hat{Z}}$. Hence
\[
\pi_{\ast}\Omega_{\hat{Z}}=\pi_{\ast}\pi^{\ast}\Omega_Z = \Omega_Z \otimes \pi_{\ast}\mathcal{O}_{\hat{Z}}
\]
Now from~(\ref{eq-ob3-0}) there is a natural surjective map
\begin{equation}\label{eq-ob3-2}
\pi_{\ast}\Omega_{\hat{Z}}=\Omega_Z \otimes \pi_{\ast}\mathcal{O}_{\hat{Z}} \la \Omega_Z\otimes L \la 0
\end{equation}
Now composing the maps~(\ref{eq-ob3-1}) and~(\ref{eq-ob3-2}) we obtain a surjective map
\begin{equation}\label{eq-ob3-3}
\delta_1 \colon \pi_{\ast} \Omega_{\tilde{X}} \la \Omega_Z \otimes L
\end{equation}
Thus we have defined the maps $\delta_0$ and $\delta_1$. To conclude it remains to show that the sequence
\begin{equation}\label{eq-ob3-4}
\Omega_X \stackrel{\delta_0}{\la} \pi_{\ast}\Omega_{\tilde{X}} \stackrel{\delta_1}{\la} \Omega_Z \otimes L
\end{equation}
is exact. The naturality of the definition of the maps $\delta_0$, $\delta_1$ allows us to check the exactness locally. So we may assume that $X=\mathrm{Spec} R$, where $R = \mathbb{C}[x_1,\ldots,x_n]/(x_1x_2)$. Then $Z=\mathrm{Spec} B$, $B=R/(x_1,x_2)$, $\tilde{X}=X_1 \coprod X_2 = \mathrm{Spec}(R_1 \times R_2)$, where $R_i=\mathbb{C}[x_1,\ldots , x_n]/(x_i)$, $i=1,2$, and $\tilde{Z}=Z_1 \coprod Z_2= \mathrm{Spec} (B_1 \times B_2)$, where $B_1 =R_1/(x_2)$ and $B_2=R_2/(x_1)$. The map $\pi \colon \tilde{X} \la X$ is induced by the map
\[
p \colon R \la R_1 \times R_2
\]
such that $p(\overline{f}(x_1,\ldots ,x_n)=(\overline{f}(0,x_2,\ldots,x_n),\overline{f}(x_1,0,x_3,\ldots,x_n))$. In turn, the map $\delta_0 \colon \Omega_X \la \pi_{\ast}\Omega_{\tilde{X}}$ is induced by the map
\[
\sigma_0 \colon \Omega_R \la \Omega_{R_1\times R_2}=\Omega_{R_1} \times \Omega_{R_2}
\]
defined by $\sigma(d\overline{f}(x_1,x_2,\ldots,x_n))=(d(\overline{f}(0,x_2,\ldots,x_n)),d(\overline{f}(x_1,0,x_3,\ldots,x_n))$, where by $\overline{f}$ we define the image of any polynomial of $\mathbb{C}[x_1,\ldots,x_n]$ to any of the rings defined earlier.

Next we describe $\delta_1$. The map $\beta \colon \mathcal{O}_Z \la \pi_{\ast} \mathcal{O}_{\hat{Z}}$ is induced from the map
\[
p_1 \colon B \la B_1 \times B_2
\]
given by $p_1(\overline{f}(x_1,\ldots ,x_n)=(\overline{f}(0,x_2,\ldots,x_n),\overline{f}(x_1,0,x_3,\ldots,x_n))$. Hence the exact sequence~(\ref{eq-ob3-0}) is induced by
\[
0 \la B \stackrel{p_1}{\la}B_1 \times B_2 \stackrel{p_2}{\la} B \la 0
\]
where $p_2$ is given by $p_2$ is given by $p_2(\overline{f}(x_1,\ldots,x_n),\overline{g}(x_1,\ldots,x_n))=\overline{f}(x_1,\ldots,x_n)-\overline{g}(x_1,\ldots,x_n)$. It is easily seen that the above sequence is an exact sequence of $B$-modules. Hence~\(\ref{eq-ob3-4})$ is induced from
\begin{equation}\label{eq-ob3-5}
\Omega_R \stackrel{\beta}{\la} \Omega_{R_1}\times \Omega_{R_2} \stackrel{\gamma}{\la} \Omega_B
\end{equation}
where
\[
\beta(d(\overline{f}(x_1,\ldots,x_n))=(d(\overline{f}(0,x_2,\ldots,x_n)) , d(\overline{f}(x_1,0,x_3,\ldots,x_n))
\]
and
\[
\gamma(d(\overline{f}(x_1,x_2,\ldots,x_n)) , d(\overline{f}(x_1,x_2,x_3,\ldots,x_n))=d(\overline{f}(x_1,x_2,\ldots,x_n)) - d(\overline{f}(x_1,x_2,x_3,\ldots,x_n))
\]
It is now straightforward to check that~(\ref{eq-ob3-5}) is exact, and therefore so is~(\ref{eq-ob3-4}) is exact too, as claimed.

\textbf{Step 2.} We now do the general case. So, let $X$ be any scheme with normal crossing singularities, $Z$ its singular locus and $\pi \colon \tilde{X} \la X$, $\nu \colon \tilde{Z} \la Z$ the normalizations of $X$ and $Z$ respectively. Following the proof of Theorem~\ref{2}, there is a sequence of morphisms
\[
X^{\prime} \stackrel{g}{\la} X_m \stackrel{f_m}{\la} \cdots \stackrel{f_2}{\la} X_1 \stackrel{f_1}{\la} X
\]
such that $f_i$ are birational for all $1\leq i \leq m$, $X_m$ has only normal crossing points of multiplicity 2 and $g$ is its normalization. Let $Z_m \subset X_m$ the birational transform of $Z$ in $X_m$. Then it is smooth and equal to the singular locus of $X_m$. Moreover, it is birational over $Z$. Let $f=f_1 \circ \cdots \circ f_m$. Then there are birational morphisms $h \colon X^{\prime} \la \tilde{X}$ and $\alpha \colon Z_m \la \tilde{Z}$ such that $\pi \circ h =f\circ g$ and $f=\nu \circ \alpha$. Let $Z^{\prime}=g^{-1}Z_m$. Then from step 1 there is an exact sequence
\begin{equation}\label{eq-ob3-6}
0 \la \Omega_{X_m}/\tau_{X_m} \la g_{\ast}\Omega_{X^{\prime}} \la \Omega_{Z_m} \otimes L_m \la 0
\end{equation}
where $L_m$ is a 2-torsion invertible sheaf on $Z_m$. Applying $f_{\ast}$ we get the exact sequence
\begin{equation}\label{eq-ob3-7}
0 \la f_{\ast}(\Omega_{X_m}/\tau_{X_m}) \la f_{\ast}g_{\ast}\Omega_{X^{\prime}} \la f_{\ast}(\Omega_{Z_m} \otimes L_m)
\end{equation}
\textbf{Claim:} \[
f_{\ast}(\Omega_{X_m}/\tau_{X_m}) = \Omega_{X}/\tau_{X}
\]
There are natural maps
\begin{gather}
\Omega_{\tilde{X}}\la h_{\ast} \Omega_{X^{\prime}}\\
\Omega_{Z_k} \la \alpha_{\ast} \Omega_{\tilde{Z}}
\end{gather}
These maps are isomorphisms since $h$ and $\alpha$ are birational and $\tilde{Z}$, $Z_m$, $X^{\prime}$, $\tilde{X}$ are smooth. The torsion sheaves at question are defined by the exact sequences
\begin{gather*}
0 \la \tau_X \la \Omega_X \la \pi_{\ast}\Omega_{\tilde{X}}\\
0 \la \tau_{X_m} \la \Omega_{X_m} \la g_{\ast}\Omega_{X^{\prime}}
\end{gather*}
Applying $f_{\ast}$ in the second we obtain a commutative diagram
\[
\xymatrix{
0 \ar[r] & \tau_X \ar[r] & \Omega_X \ar[r]\ar[d] & \pi_{\ast}\Omega_{\tilde{X}} \ar@{=}[d] \\
0 \ar[r] & f_{\ast}\tau_{X_m} \ar[r]&     f_{\ast}\Omega_{X_m} \ar[r] &  f_{\ast}g_{\ast}\Omega_{X^{\prime}}=\pi_{\ast}h_{\ast}\Omega_{X^{\prime}}=\pi_{\ast}\Omega_{\tilde{X}}
}
\]
Hence there is a map $\tau_X \la f_{\ast}\tau_{X_m}$ which induces a natural map $\Omega_{X}/\tau_{X} \la f_{\ast}(\Omega_{X_m}/\tau_{X_m})$. Since $ \Omega_{X}/\tau_{X}$ is reflexive by Proposition~\ref{ob1} and $f$ is birational, the map is an isomorphism as claimed. Moreover, since $\alpha$ is birational and $L_m$ is 2-torsion, $L_m=\alpha^{\ast}L$, where $L$ is a 2-torsion invertible sheaf on $\tilde{Z}$.  Now
\[
f_{\ast}g_{\ast}\Omega_{X^{\prime}}=\pi_{\ast}h_{\ast}\Omega_{X^{\prime}}=\pi_{\ast}\Omega_{\tilde{X}}
\]
and
\[
f_{\ast}(\Omega_{Z_m}\otimes L_m) =\nu_{\ast}\alpha_{\ast}(\Omega_{Z_m}\otimes \alpha^{\ast}L) =\nu_{\ast}(\Omega_{\tilde{Z}}\otimes L)
\]
Therefore,~(\ref{eq-ob3-7}) becomes
\[
0 \la \Omega_X/\tau_X \la \pi_{\ast}\Omega_{\tilde{X}} \la \nu_{\ast}(\Omega_{\tilde{Z}}\otimes L)
\]
as was to be proven.
\end{proof}

\begin{theorem}\label{ob4}
Let $X$ be a scheme with normal crossing singularities such that $\dim X \leq 3$ and $\omega_X^{-1}$ is ample. Then $H^2(T_X)=0$.
\end{theorem}

\begin{proof}
We only do the 3-fold case. The surface case is exactly similar and is omitted. So, let $X$ be a Fano 3-fold with normal crossing singularities and let $Z \subset X$ be its singular locus. Then by Corollary~\ref{ob2},
\[
H^2(T_X)=H^1((\Omega_X/\tau_X)\otimes \omega_X)
\]
Then by Proposition~\ref{ob3}, there is an exact sequence
\[
0 \la \Omega_X/\tau_X \stackrel{\delta_1}{\la} \pi_{\ast}\Omega_{\tilde{X}} \stackrel{\delta_2}{\la} \nu_{\ast}(\Omega_{\tilde{Z}}\otimes L)
\]
where $\pi \colon \tilde{X} \la X$, $\nu \colon \tilde{Z} \la Z$ are the normalizations of $X$ and $Z$, respectively and $L$ is a 2-torsion invertible sheaf on $\tilde{Z}$. Tensoring it with $\omega_X$ and taking into consideration that $\pi_{\ast}(\Omega_{\tilde{X}} \otimes \pi^{\ast}\omega_X)=\pi_{\ast} \Omega_{\tilde{X}} \otimes \omega_X$, we get the exact sequence
\[
0 \la (\Omega_X/\tau_X)\otimes \omega_X  \stackrel{\delta_1}{\la} \pi_{\ast}(\Omega_{\tilde{X}}\otimes \pi^{\ast}\omega_X) \stackrel{\delta_2}{\la} \nu_{\ast}(\Omega_{\tilde{Z}}\otimes L \otimes  \nu^{\ast}\omega_X)
\]
Let $N=\mathrm{Coker}(\delta_1)$. Then the above sequence splits into
\begin{gather*}
0 \la (\Omega_X/\tau_X)\otimes \omega_X  \stackrel{\delta_1}{\la} \pi_{\ast}(\Omega_{\tilde{X}}\otimes \pi^{\ast}\omega_X) \stackrel{\delta_2}{\la} N\la 0\\
0 \la N \la \nu_{\ast}(\Omega_{\tilde{Z}}\otimes L \otimes \nu^{\ast}\omega_X)
\end{gather*}
Therefore we an exact sequence
\begin{equation}\label{eq-ob4-1}
H^0(N) \la H^1((\Omega_X/\tau_X)\otimes \omega_X) \la H^1(\pi_{\ast}(\Omega_{\tilde{X}}\otimes \pi^{\ast}\omega_X))
\end{equation}
Now since $\pi$ and $\nu$ are finite, it follows that $(\pi^{\ast}\omega_X)^{-1}$ and $(\nu ^{\ast}\omega_X)^{-1}$ are ample. Moreover, since $L$ is 2-torsion and invertible, $(L^{-1} \otimes\nu ^{\ast}\omega_X)^{-1}$ is ample too. Therefore, and by using the Kodaira-Nakano vanishing theorem,
\[
H^0(N)\subset H^0(\nu_{\ast}(\Omega_{\tilde{Z}}\otimes L \otimes  \nu^{\ast}\omega_X))=H^0(\Omega_{\tilde{Z}}\otimes L \otimes \nu^{\ast}\omega_X)=0
\]
and
\[
H^1(\pi_{\ast}(\Omega_{\tilde{X}}\otimes \pi^{\ast}\omega_X))=H^1(\Omega_{\tilde{X}}\otimes \pi^{\ast}\omega_X)=0
\]
Hence from~(\ref{eq-ob4-1}) it follows that
\[
H^2(T_X)=H^1((\Omega_X/\tau_X)\otimes \omega_X) =0
\]
as claimed.
\end{proof}

Unfortunately, in general I cannot say much about the other obstruction space, namely $H^1(T^1(X))$. However, since $T^1(X)$ is a line bundle on the singular locus $Z$ of $X$, it is much easier handled than $H^2(T_X)$ and will vanish if we impose certain positivity requirement on $T^1(X)$.

The case when $X$ has only double points exhibits much better behaviour and it deserves special consideration. The difference between this and the general case is that the exact sequence~\ref{ob-seq} in Theorem~\ref{ob3} is right exact too which makes it possible to calculate obstructions in all dimensions.

\subsection{The double points case.}

\begin{theorem}\label{db-ob}
Let $X$ be a scheme with only double point normal crossing singularities such that $\omega_X^{-1}$ is ample. Then
\begin{enumerate}
\item \[
H^2(T_X)=0
\]
\item Suppose that $T^1(X)$ is finitely generated by its global sections. Then\[
H^1(T^1(X))=0
\]
\end{enumerate}
\end{theorem}
\begin{corollary}
Let $X$ be a scheme with only double point normal crossing singularities such that $\omega_X^{-1}$ is ample and such that $T^1(X)$ is finitely generated by its global sections. Then $\mathrm{Def}(X)$ is smooth.
\end{corollary}
\begin{question}
Is $Def(X)$ smooth for any Fano variety with normal crossing singularities? If this is true then $X$ is smoothable if and only if $T^1(X)$ is finitely generated by global sections and hence this is a very natural condition to impose.
\end{question}
\begin{remark}
In general, $H^1(T^1(X))$ will not vanish. However if $X$ is smoothable, then $T^1(X)$ must have some positivity properties and the one stated is the most natural one.
\end{remark}
\begin{proof}[Proof of Theorem~\ref{db-ob}]
The first part is proved in excactly the same way as Theorem~\ref{ob4}. The crucial part in the double point case is that the exact sequence~(\ref{ob-seq}) is right exact too.

In order to show the second part we will show that the singular locus $Z$ of $X$ is a smooth Fano variety of dimension $\dim X-1$. The Fano part is the only part to be shown. Let $\pi \colon \tilde{X} \la X$ be the normalization and $\tilde{Z}=\pi^{-1}Z$. Then $\tilde{Z} \la Z$ is \'etale. By subadjunction we get that
\[
\pi^{\ast}\omega_X=\omega_{\tilde{X}} \otimes \mathcal{O}_{\tilde{X}}(\tilde{Z})
\]
Therefore,
\[
\omega_{\tilde{Z}}=\omega_{\tilde{X}}\otimes \mathcal{O}_{\tilde{X}}(\tilde{Z})\otimes \mathcal{O}_{\tilde{Z}}
\]
Hence $\omega_{\tilde{Z}}^{-1}$ is ample. But since $\tilde{Z} \la Z$ is \'etale, it follows that $\pi^{\ast}\omega_Z=\omega_{\tilde{Z}}$. Therefore $\omega_Z^{-1}$ is ample too and hence $Z$ is Fano as claimed.  Now
\[
H^1(T^1(X))=H^1(\omega_Z \otimes(T^1(X)\otimes \omega_Z^{-1}))=0
\]
by the Kawamata-Viehweg vanishing theorem since if $T^1(X)$ is finitely generated by global sections, then $T^1(X) \otimes \omega_Z^{-1}$ is ample too.

\end{proof}

\section{Smoothings of Fanos}
In this section we obtain criteria for a Fano variety $X$ with normal crossing singularities to be smoothable. 
\begin{proposition}\label{sm1}
Let $X$ be a reduced projective scheme with hypersurface singularities. Then
\begin{enumerate}
\item If $X$ is smoothable by a flat deformation $\mathcal{X} \la \Delta$ such that $\mathcal{X}$ is smooth, then $T^1(X)=\mathcal{O}_Z$, where $Z$ is the singular locus of $X$.
\item Suppose that $T^1(X)$ is finitely generated by its global sections and that $H^2(T_X)=H^1(T^1(X))=0$. Then $X$ is smoothable. Moreover, if $Def(X)$ is smooth then the converse is also true.
\end{enumerate}
\end{proposition}

\begin{proof}
The second part is~\cite[Theorem 12.5]{Tzi09a}. We proceed to show the first part. Let $f \colon \mathcal{X} \la \Delta$ be a smoothing of $X$ such that $\mathcal{X}$ is smooth, where $\Delta=\mathrm{Spec}(R)$, $(R,m_R)$ is a discrete valuation ring. Let $T^1(\mathcal{X}/\Delta)=\mathcal{E}xt^1_{\mathcal{X}}(\Omega_{\mathcal{X}/\Delta},\mathcal{O}_{\mathcal{X}})$ be Schlessinger's relative $T^1$ sheaf. Then dualizing  the exact sequence
\[
0 \la f^{\ast}\omega_{\Delta} = \mathcal{O}_{\mathcal{X}} \la \Omega_{\mathcal{X}} \la \Omega_{\mathcal{X}/\Delta} \la 0
\]
we get the exact sequence
\[
\cdots \la \mathcal{O}_{\mathcal{X}} \la T^1(\mathcal{X}/\Delta) \la \mathcal{E}xt^1_{\mathcal{X}}(\Omega_{\mathcal{X}},\mathcal{O}_{\mathcal{X}}) =0
\]
Now restricting to the special fiber and taking into consideration that $\mathcal{X} \otimes_R R/m_r \cong X$ and that $T^1(\mathcal{X}/\Delta)\otimes_R R/m_R=T^1(X)$~\cite[Lemma 7.7]{Tzi09a}, we get that there is a surjection $\mathcal{O}_X \la T^1(X)$. Moreover, $T^1(X)$ is a line bundle on the singular locus $Z$ of $X$. Hence, restricting on $Z$ it follows that $T^1(X) \cong \mathcal{O}_Z$, as claimed.
\end{proof}

\begin{remark}
The condition $T^1(X)=\mathcal{O}_Z$ is equivalent to Friedman's d-semistability condition in the case of reducible simple normal crossing schemes~\cite{Fr83}. One of the natural questions raised by Friedman is whether this condition is sufficient for a simple normal crossing variety to be smoothable. He showed that in the case of $K3$ surfaces this is true but Persson and Pinkham have shown that this is not true in general~\cite{PiPe83}. However this is true in the case of normal crossing (not necessarily reducible) Fano schemes with at most double points, as shown by Theorem~\ref{sm3} below.
\end{remark}

\begin{theorem}\label{sm2}
Let $X$ be a scheme with normal crossing singularities of dimension at most 3. Suppose that $T^1(X)$ is generated by global sections and that $H^1(T^1(X))=0$. Then $X$ is smoothable.
\end{theorem}

\begin{proof}
The theorem follows directly from Theorem~\ref{ob4} and Proposition~\ref{sm1}.2.
\end{proof}

\begin{theorem}\label{sm3}
Let $X$ be a scheme, of any dimension, with double point normal crossing singularities. Then
\begin{enumerate}
\item If $T^1(X)$ is finitely generated by its global sections, then $X$ is smoothable.
\item $X$ is smoothable by a flat deformation $\mathcal{X} \la \Delta$ with $\mathcal{X}$ smooth, if and only if $T^1(X)=\mathcal{O}_Z$, where $Z$ is the singular locus of $X$.
\end{enumerate}
\end{theorem}

\begin{proof}
The theorem follows immediately from Theorem~\ref{db-ob} and Proposition~\ref{sm1}.
\end{proof} 
\section{Simple normal crossing Fanos.}
As mentioned in the introduction, the case of simple normal crossing Fano varieties deserves to be treated separately. So let $X$ be a simple normal crossing Fano scheme of dimension $n$ and let $X_i$, $1 \leq i \leq k$, be its irreducible components, which by definition are smooth. In this section we study when $X$ is smoothable and in particular when it is smoothable by a flat morphism $f \colon \mathcal{X} \la \Delta$ such that $\mathcal{X}$ is smooth. The main difference between this and the general case is that now we can use the general theory on simple normal crossing varieties developed by Friedman~\cite{Fr83} and the theory of logarithmic deformations developed by Kawamata and Namikawa in the complex analytic case~\cite{KawNa94} and K. Kato, F. Kato in the algebraic case~\cite{Ka88},~\cite{Kat96}.

First we study the problem of when $X$ is smoothable and then when it is smoothable with a smooth total space. Following the Friedman~\cite{Fr83} we define
\[
\alpha_p \colon X^{[p]}=\coprod_{i_0 < \cdots < i_p} (X_{i_0}\cap \cdots \cap X_{i_p} ) \la X
\]
Note that $X^{[0]}$ is the normalization $\tilde{X}$ of $X$ and $X^{[1]}$ is the normalization $\tilde{Z}$ of the singular locus $Z$ of $X$. Then $X^{i]}$ is smooth of dimension $n-i$, for all $i \leq n$ and there is an exact sequence~\cite{Fr83}
\begin{equation}\label{sec5-eq1}
0 \la \tau_X \la \Omega_X \la (\alpha_0)_{\ast} \Omega_{X^[0]} \la (\alpha_1)_{\ast} \Omega_{X^[1]} \la \cdots
\end{equation}
which extends the exact sequence~(\ref{ob-seq}). Then arguing in exactly the same way as in the proof of Theorem~\ref{ob4} and Theorem~\ref{sm2} we find that
\begin{theorem}\label{snc-ob}
Let $X$ be a simple normal crossing Fano scheme. Then $H^2(T_X)=0$. Hence if $T^1(X)$ is finitely generated by global sections and $H^1(T^1(X))=0$, then $X$ is smoothable.
\end{theorem}

Next we study the problem of when $X$ admits a smoothing with a smooth total space. In what follows we use the theory of logarithmic deformations. For details and basic properties we refer the reader to~\cite{KawNa94} and~\cite{Ka88},~\cite{Kat96}.

Suppose that $f \colon \mathcal{X} \la \Delta$ is a smoothing such that $\mathcal{X}$ is smooth. Then by Proposition~\ref{sm1}, $T^1(X) = \mathcal{O}_Z$, where $Z$ is the singular locus of $X$. We will study to what extend the converse is also true. According to Friedman~\cite{Fr83}, if $T^1(X)=\mathcal{O}_Z$, then $X$ is $d$-semistable. Therefore $X$ admits a logarithmic structure~\cite{KawNa94}~\cite{Kat96}. Then there is a locally free sheaf of logarithmic deformations $\Omega_{X/\mathbb{C}}(\mathrm{log})$ such that $H^2(T_{X/\mathbb{C}}(\mathrm{log}))$ is an obstruction space for logarithmic deformations~\cite[Theorem 2.2]{KawNa94}, where $T_{X/\mathbb{C}}(\mathrm{log})=\mathcal{H}om_X(\Omega_{X/\mathbb{C}}(\mathrm{log}),\mathcal{O}_X)$.

The next result allows us to make calculations involving $\Omega_{X/\mathbb{C}}(\mathrm{log})$.

\begin{proposition}[~\cite{Fr83}]\label{log-ex}
The following sequence is exact
\[
0 \la \Omega_X/\tau_X  \la \Omega_{X/\mathbb{C}}(\mathrm{log}) \la (\alpha_1)_{\ast}\mathcal{O}_{X^{[1]}} \la (\alpha_2)_{\ast}\mathcal{O}_{X^{[2]}} \la \cdots
\]
\end{proposition}

We are now able to calculate the obstruction to logarithmic deformations.
\begin{theorem}\label{snc-sm}
Let $X$ be a simple normal crossing Fano scheme. Then
\[
H^2(T_{X/\mathbb{C}}(\mathrm{log}))=0.
\]
\end{theorem}
\begin{proof}
The proof is similar to the proof of Theorem~\ref{ob4} and so we simply describe the main steps.
By Serre duality $H^2(T_{X/\mathbb{C}}(\mathrm{log}))=H^{n-2}(\Omega_{X/\mathbb{C}}(\mathrm{log})\otimes \omega_X)$, where $n=\dim X$. Now from~(\ref{sec5-eq1}) and by using the Kodaira-Nakano vanishing theorem we get that $H^{n-2}((\Omega_X/\tau_X) \otimes \omega_X)=0$. Now from Proposition~\ref{log-ex} and by using Kodaira vanishing we find that $
H^2(T_{X/\mathbb{C}}(\mathrm{log}))=H^{n-2}(\Omega_{X/\mathbb{C}}(\mathrm{log}) \otimes \omega_X)=0$.

\end{proof}
We now use~\cite[Corollary 2.4]{KawNa94} to conclude that any d-semistable simple normal crossing Fano scheme is smoothable.
\begin{theorem}\label{snc-sm1}
Let $X$ be a simple normal crossing scheme. Then $X$ has a smoothing $f \colon \mathcal{X}\la \Delta$ such that $\mathcal{X}$ is smooth, if and only if $T^1(X)\cong \mathcal{O}_Z$, where $Z$ is the singular locus of $X$.
\end{theorem}

\section{Examples.}
In this section we construct one example of a smoothable and one of a non-smoothable normal crossing Fano 3-fold.

\begin{example}\label{ex1} Let $P \in Y \subset \mathbb{P}^4$ be a quadric surface with one ordinary double point locally analytically isomorphic to $(xy-zw=0)\subset \mathbb{C}^4$. Let $f \colon X \la Y$ be the blow up of $P \in Y$. Then $X$ is smooth and the $f$-exceptional divisor $E$ is isomorphic to $\mathbb{P}^1 \times \mathbb{P}^1$. Moreover, $-K_X-E$ is ample and $\mathcal{N}_{E/X}=\mathcal{O}_E(-1,-1)$.

Next we construct an embedding $E \subset X^{\prime}$ of $E$ into a smooth scheme $X^{\prime}$ such that, $\mathcal{N}_{E/X^{\prime}}=\mathcal{O}_{E}(1,1)$ and $-K_{X^{\prime}}-E$ is ample. Let $Z \subset \mathbb{P}^3$ be a smooth quadric surface. Then $\mathcal{N}_{Z/\mathbb{P}^3}=\mathcal{O}_Z(2,2)$. Let $\pi \colon X^{\prime} \la \mathbb{P}^3$ be the cyclic double cover of $\mathbb{P}^3$ ramified over $Z$. This is defined by the line bundle $\mathcal{L}=\mathcal{O}_{\mathbb{P}^3}(1)$ and the section $s$ of $\mathcal{L}^{\otimes 2}$ that corresponds to $Z$. Let $E=(\pi^{-1}(Z))_{\mathrm{red}} \cong Z$. Then $\pi^{\ast}Z = 2E$ and $\omega_{X}=\pi^{\ast}(\omega_{\mathbb{P}^3}\otimes \mathcal{L})$. Let $l^{\prime} \subset E$ be one of the rulings and $l = \pi_{\ast}(l^{\prime})$. then
\[
l^{\prime} \cdot E = 1/2 (l^{\prime} \cdot \pi^{\ast} Z) = 1/2( l \cdot Z) = 1
\]
and hence $\mathcal{N}_{E/\mathcal{X}^{\prime}}=\mathcal{O}_{E}(1,1)$.  Now let $Y$ be the scheme obtained by glueing $X$ and $X^{\prime}$ along $E$. This is a normal crossing Fano 3-fold with only double points. Then by Theorem~\ref{1}, $T^1(Y)= \mathcal{N}_{E/X} \otimes \mathcal{N}_{E/X^{\prime}} = \mathcal{O}_{E}$. Therefore, by Theorem~\ref{sm3}, $Y$ is smoothable.
\end{example}

\begin{example}\label{ex2} Let $E \subset X$ be as in example 1. Then let $Y$ be obtained by glueing two copies of $X$ along $E$. Then by Theorem~\ref{1},
\[
T^1(Y)=\mathcal{N}_{E/X} \otimes \mathcal{N}_{E/X} = \mathcal{O}_E(-2,-2)
\]
and hence $H^0(T^1(Y))=0$. Hence $Y$ is not smoothable~\cite[Theorem 12.3]{Tzi09}. In fact, every deformation of $Y$ is locally trivial.
\end{example}

\end{document}